\documentclass[12pt]{article}
\usepackage{natbib}
\usepackage{color}
\usepackage{latexsym}
\usepackage{amsthm}
\usepackage{multirow}
\usepackage{float}
\usepackage{wrapfig}
\usepackage[latin1]{inputenc}
\usepackage[T1]{fontenc}
\usepackage{amsmath,amsfonts,amssymb,bm,dsfont}
\usepackage{pslatex}
\usepackage{graphicx,color}
\usepackage{fancyhdr}
\usepackage{multicol}
\usepackage{thumbpdf,lmodern}

\makeatletter
\renewcommand\@biblabel[1]{}

\makeatother
\usepackage[displaymath]{lineno}
\usepackage[latin1]{inputenc}
\usepackage[T1]{fontenc}
\usepackage{graphicx}
\usepackage{times}
\usepackage{amsmath,amsfonts,amssymb,bm}
\usepackage[english]{babel}
\usepackage{amsbsy}

\usepackage{subfigure}
\usepackage{multirow}
\usepackage{multicol}
\theoremstyle{definition}

 \newcommand*\patchAmsMathEnvironmentForLineno[1]{%
  \expandafter\let\csname old#1\expandafter\endcsname\csname #1\endcsname
  \expandafter\let\csname oldend#1\expandafter\endcsname\csname end#1\endcsname
  \renewenvironment{#1}%
     {\linenomath\csname old#1\endcsname}%
     {\csname oldend#1\endcsname\endlinenomath}}%
\newcommand*\patchBothAmsMathEnvironmentsForLineno[1]{%
  \patchAmsMathEnvironmentForLineno{#1}%
  \patchAmsMathEnvironmentForLineno{#1*}}%
\AtBeginDocument{%
\patchBothAmsMathEnvironmentsForLineno{equation}%
\patchBothAmsMathEnvironmentsForLineno{align}%
\patchBothAmsMathEnvironmentsForLineno{flalign}%
\patchBothAmsMathEnvironmentsForLineno{alignat}%
\patchBothAmsMathEnvironmentsForLineno{gather}%
\patchBothAmsMathEnvironmentsForLineno{multline}%
}

\begin{document}

\begin{center}

\Large {\bf Estimation of entropy measures for categorical variables with spatial correlation}\par

\bigskip

\normalsize{Linda Altieri, Daniela Cocchi, Giulia Roli}\\ \small{Department of Statistical Sciences, University of Bologna, via Belle Arti, 41, 40126, Bologna, Italy.} \par

\bigskip

%

\end{center}

\begin{quotation}
\noindent {\it Abstract:}
Entropy is a measure of heterogeneity widely used in applied sciences, often when data are collected over space. Recently, a number of approaches has been proposed to include spatial information in entropy. The aim of entropy is to synthesize the observed data in a single, interpretable number. In other studies the objective is, instead, to use data for entropy estimation; several proposals can be found in the literature, which basically are corrections of the estimator based on substituting the involved probabilities with proportions. In this case, independence is assumed and spatial correlation is not considered.
We propose a path for spatial entropy estimation: instead of correcting the global entropy estimator, we focus on improving the estimation of its components, i.e. the probabilities, in order to account for spatial effects. Once probabilities are suitably evaluated, estimating entropy is straightforward since it is a deterministic function of the distribution. Following a Bayesian approach, we derive the posterior probabilities of a multinomial distribution for categorical variables, accounting for spatial correlation. A posterior distribution for entropy can be obtained, which may be synthesized as wished and displayed as an entropy surface for the area under study.
\bigskip
\end{quotation}

\begin{quotation}
\bigskip
\end{quotation}

\section{Introduction}

Shannon's entropy is a successful measure in many fields, as it is able to synthesize several concepts in a single number: entropy, information heterogeneity, surprise, contagion. The entropy of a categorical variable $X$ with $I<\infty$ outcomes is
\begin{equation}
H(X)=\sum_{i=1}^I \log(p(x_i)) \log\left(\frac{1}{p(x_i)}\right),
\label{eq:shan}
\end{equation}
where $p(x_i)$ is the probability associated to the $i$-th outcome \citep{coverthomas}. The flexibility of such index and its ability to describe any kind of data, including categorical variables, motivates its diffusion across applied fields such as geography, ecology, biology and ladscape studies. Often, such disciplines deal with spatial data, and the inclusion of spatial information in entropy measures has been the target of intensive research (see, e.g., \citeauthor{batty76}, \citeyear{batty76}, \citeauthor{leibovici09}, \citeyear{leibovici09}, \citeauthor{leibovici14}, \citeyear{leibovici14}). 
In several case studies, the interest lies in describing and synthesizing what is observed. This is usually not a simple task: large amounts of data require advanced computational tools, qualitative variables have limited possibilities, and, when data are georeferenced, spatial correlation should be accounted for. When it comes to measuring the entropy of spatial data, we suggest an approach proposed in \cite{nostro}, which allows to decompose entropy into a term quantifying the spatial information, and a second term quantifying the residual heterogeneity. 

In other cases, though, the aim lies in estimating the entropy of a phenomenon, i.e. in making inference rather than description. Under this perspective, a stochastic process is assumed to generate the data according to an unknown probability function and an unknown entropy. One realization of the process is observed and employed to estimate such entropy. The standard approach relies on the so-called 'plug-in' estimator, presented in \cite{paninski}, which substitutes probabilities with observed relative frequencies in the computation of entropy:
\begin{equation}
\widehat{H}_p(X)=\sum_{i=1}^I \log(\widehat{p}(x_i)) \log\left(\frac{1}{\widehat{p}(x_i)}\right),
\label{eq:plugin}
\end{equation}
where $\widehat{p}(x_i)=n_i/n$ is the relative amount of observations of category $i$ over $n$ data. It is the non-parametric as well as the maximum likelihood estimator \citep{paninski}, and performs well when $I<\infty$ is known \citep{antos}. For unknown or infinite $I$, (\ref{eq:plugin}) is known to be biased; the most popular proposals at this regard consist of corrections of the plug-in estimator: see, for example, the Miller-Madow \citep{miller} and the jack-knifed corrections \citep{efron}. Recently, \cite{zhang12} proposed a non-parametric solution with faster decaying bias and upper limit for the variance when $I=\infty$. Under a Bayesian framework, the most widely known proposal is the NSB estimator \citep{nemenman}, improved by \cite{archer} as regards the prior distribution. In all the above mentioned works, independence among realizations is assumed. Other approaches, linked to machine learning methods, directly estimate entropy relying on the availability of huge amounts of data \citep{hausser}.

Two main limits concern entropy estimation. Firstly, the above mentioned proposals referring to (\ref{eq:plugin}) focus on correcting or improving the performance of (\ref{eq:plugin}), while different perspectives might be considered. Secondly, no study is available about estimating entropy for variables presenting spatial association: the assumption of independence is never relaxed, while spatial entropy studies do not consider inference. 

In this paper, we take a different perspective to entropy estimation, which moves the focus from the index itself to its components. As can be seen from (\ref{eq:shan}), entropy is a deterministic function of the probability mass function (pmf) of the variable of interest. Therefore, once the pmf is properly estimated, the subsequent steps are straightforward. In the case of categorical variables following, e.g., a multinomial distribution, the crucial point is to estimate the distribution parameters. A Bayesian approach allows to derive the pmf of such distribution, and can be extended to account for spatial correlation among categories. After obtaining a posterior distribution for the parameters, the posterior distribution of entropy is also available as a transformation. Thus, a point estimator of entropy can be, e.g., the mean of the posterior distribution of the transformation; credibility intervals and other syntheses may be obtained via the standard tools of Bayesian inference. This approach can be used for non-spatial settings as well; in the spatial context, coherently with standard procedures for variables linked to areal and point data, the estimation output is a smooth spatial surface for the entropy over the area under study.

The paper is organized as follows. Section \ref{sec:theory} summarizes the methodology for Bayesian spatial regression and shows how to obtain the posterior distribution and the Bayesian estimator of entropy. Then, Section \ref{sec:data} assesses the performance of the proposed method on simulated data for different spatial configurations. Lastly, Section \ref{sec:disc} discusses the main results.

\section{Bayesian spatial entropy estimation}
\label{sec:theory}

For simplicity of presentation, we focus on the binary case. Let $X$ be a binary variable with $x_1=1$ and $x_2=2$; consider a series of $n$ realizations indexed by $u=1, \dots,n$, each carrying an outcome $x_u \in \{1,2\}$. This may be thought of as a $n$-variate variable, or alternatively as a sequence of variables $X_1, \dots X_n$, which are independent, given the distribution parameters and any effects modelling them. For a generic $X_u$, the simplest model is:
\begin{equation}
X_u \sim Ber(p_u)
\end{equation}
\begin{equation}
logit(p_u)=z_u^{'} \beta
\end{equation}
in absence of random effects, where $z_u$ are the covariates associated to the $u$-th unit. 

To the aim of including spatial correlation, consider $n$ realizations from a binary variable over the two-dimensional space, where $u=(x,y)$ identifies a unit via its spatial coordinates. Let us consider the case of realizations over a regular lattice of size $n=n_1 \times n_2$, where $u$ identifies each cell centroid. The sequence $X_1, \dots, X_n$ is now no longer independent, but spatially correlated. In order to define the extent of the correlation for grid data, the notion of neighbourhood must be introduced, linked to the assumption that occurrences at certain locations are influenced by what happens at surrounding locations, i.e. their neighbours. The simplest way of representing a neighbourhood system is via an adjacency matrix: for $n$ spatial units, $A=\{a_{uu'}\}_{u,u'=1,\dots,n}$ is a square $n\times n$ matrix such that $a_{uu'}=1$ when unit $u$ and unit $u'$ are neighbours, and $a_{uu'}=0$ otherwise; in other words, $a_{uu'}=1$ if $u' \in \mathcal{N}(u)$, the neighbourhood of area $u$, and diagonal elements are all zero by default. In the remainder of the paper, the word 'adjacent' is used accordingly to mean 'neighbouring', even when this does not correspond to a topological contact. The most common neighbourhood systems for grid data are the '4 nearest neighbours', i.e. a neighbourhood formed by the 4 pixels sharing a border along the cardinal directions, and the '12 nearest neighbours', i.e. two consequent pixels along each cardinal direction plus the four ones along the diagonals.

Auto-models provide a way of including spatial correlation: they construct a statistical model that explains a response via the response values of its neighbours. It is thus developed by combining logistic regression model with autocorrelation effects. They were initially developed for the analysis of plant competition experiments and then extended to spatial data in general. BESAG (1974) proposed to model spatial dependence among random variables directly (rather than hierarchically) and conditionally (rather than jointly). The autologistic model for spatial data with binary responses emphasises that the explanatory variables are the surrounding array variables themselves; a joint Markov random field is imposed for the binary data. A recent variant of this model, substituting Equation (4) with (5), is proposed by CARAGEA and Kaiser (2009):
\begin{equation}
logit(p_u)=z(u)'\beta + \sum_{i\in \mathcal{N}(u)} \eta_i(X_i-\mu_i)
\end{equation}
where $\eta$ parametrizes dependence on the neighbourhood and, in the simplest case, $z(u)'\beta=\beta_0$ only includes an intercept. Parameter $\mu_i=\exp(z(u)'\beta)/(1+\exp(z(u)'\beta))$ represents the expected probability of success in the situation of spatial independence.

An analogous and even more recent formulation of the model [REFERENCE], in absence of covariate information, is 
\begin{equation}
\begin{split}
& logit(p_u)=\beta_0+ \phi_u\\
& \phi \sim MVN_n(0, \Sigma)\\
& \Sigma=\left[\tau(D-\rho A)\right]^{-1}
\label{eq:model}
\end{split}
\end{equation}
where $\phi=(\phi_1, \dots, \phi_n)'$ is a spatial effect with a structured covariance matrix $\Sigma$, which depends on a precision parameter $\tau$ and a dependence parameter $\rho \in [-1, 1]$ quantifying the strength and type of the correlation between neighbouring units. $A$ is the adjacency matrix reflecting the neighbourhood structure, and $D$ is a diagonal matrix, where each element contains the row sums of $A$.

\subsection{Entropy estimation}
\label{sec:ent_est}

The estimation of the parameters for Bayesian spatial logit regression models may proceed via MCMC methods or the INLA approach. We exploit the latter \citep{rue} and obtain a posterior distribution for the parameters of the probability of success for each grid cell. A synthesis, such as the posterior mean, is chosen in order to obtain an estimate for $p_u$ over each cell. Such estimate is used for the computation of a local estimated entropy value for each pixel:
\begin{equation}
\widehat{H}(X)_u=\hat{p}_u\log\left(\frac{1}{\hat{p}_u}\right)+(1-\hat{p}_u)\log\left(\frac{1}{1-\hat{p}_u}\right).
\label{eq:ent_est}
\end{equation} 
This way, an entropy surface is obtained for estimating the process entropy, whose smoothness may be tuned by the neighbourhood choice, or by the introduction of splines for the spatial effect. Any other surface can be obtained following the same approach for different aims, e.g. for plotting the entropy standard error or the desired credibility interval extremes.

\section{Simulation study}
\label{sec:data}

\begin{figure}
\centering
     \includegraphics[width=.6\textwidth]{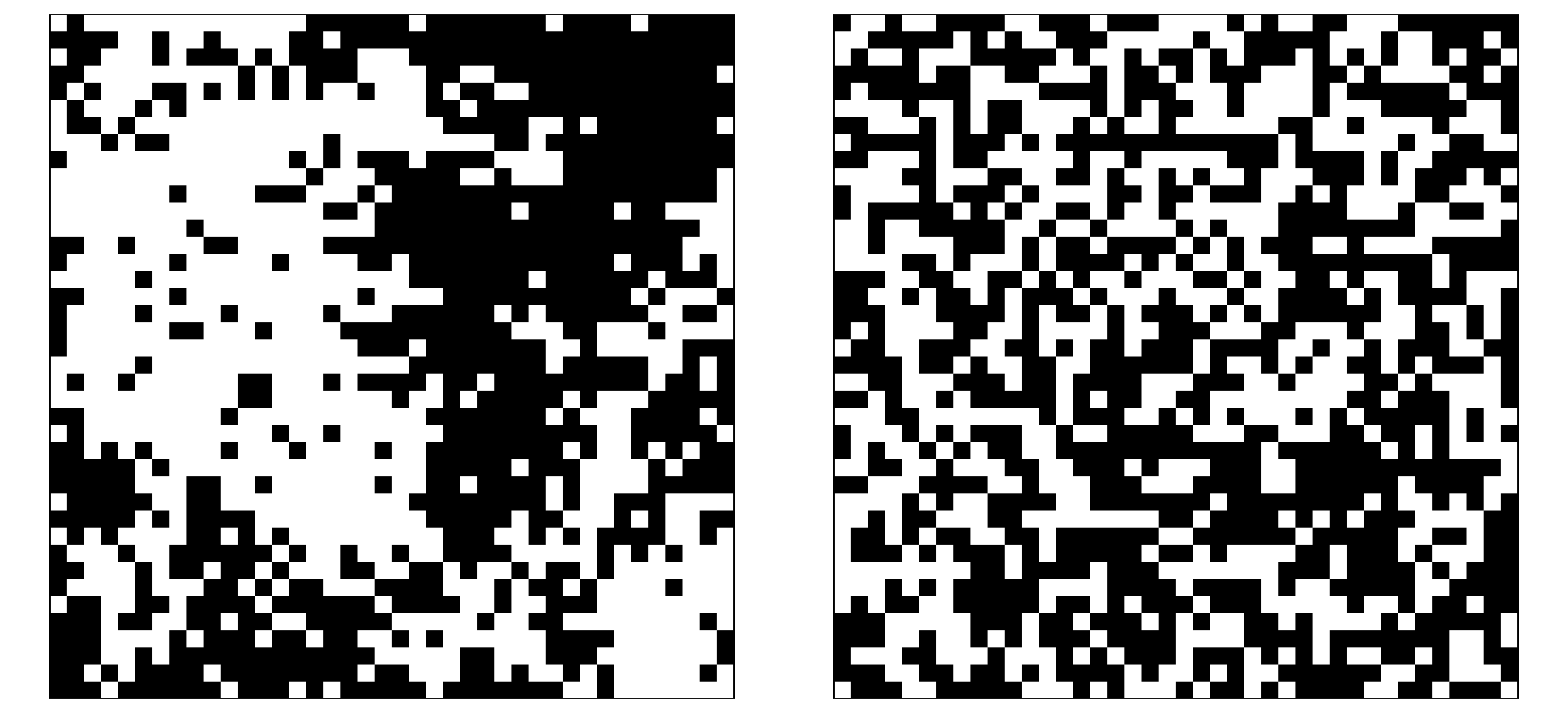}
     \caption{Clustered scenario (left) and random scenario (right) - example with $\beta_0=0.32$.}\label{fig:data_insieme}
\end{figure}

To the aim of assessing the performance of the proposed entropy estimator, we generate binary data on a $40 \times 40$ grid under two spatial configurations: clustered and random. Figure \ref{fig:data_insieme} shows an example of the generated datasets. The underlying model is model (\ref{eq:model}), with a 12 nearest neighbour structure for $A$, $\tau=0.1$ and $\rho=\{0.99, 0.0001\}$ for the two scenarios, respectively. For each scenario, 200 datasets are generated with varying values for $\beta_0$ so that the expectation of $p_u$ in a situation of independence varies between $0.1$ and $0.9$; values for $\beta_0$ differ across replicates but are constant across scenarios, so that the proportion of pixels of each type is comparable. 

\begin{figure}
\centering
     \includegraphics[width=.7\textwidth]{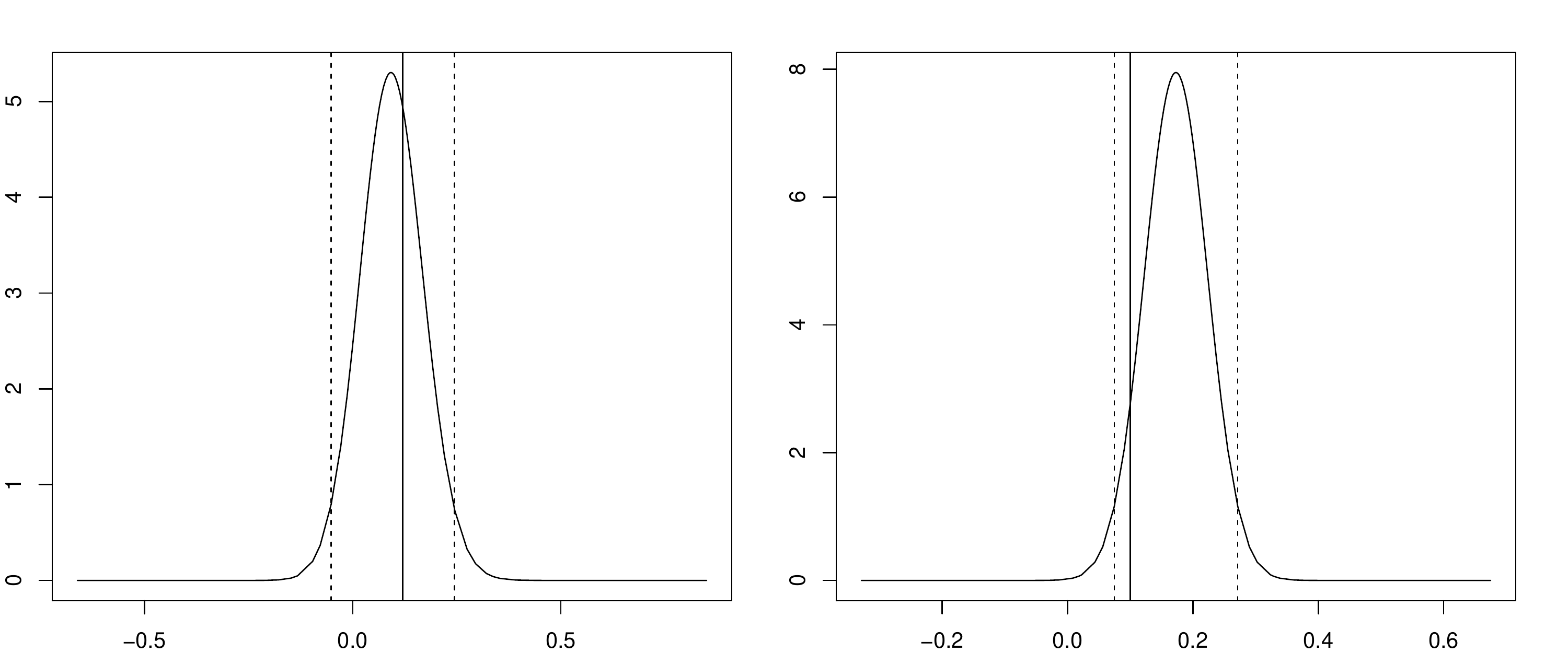}
     \caption{Posterior distributions for the model parameters, true values (vertical line) and 95\% credibility intervals (dashed vertical lines).}\label{fig:pars}
\end{figure}

Results show that fitting the model over the generated data leads to good estimates for the $p_u$s. For all scenarios, the estimated parameters are very close to the true ones, which are always included within the 95\% credibility intervals. An example is shown in Figure \ref{fig:pars}.
The proposed approach is able to produce good estimates for the probabilities of success, ensuring the goodness of the estimates for spatial entropy, which is a function of such probabilities. 

\begin{figure}
\centering
     \includegraphics[width=.7\textwidth]{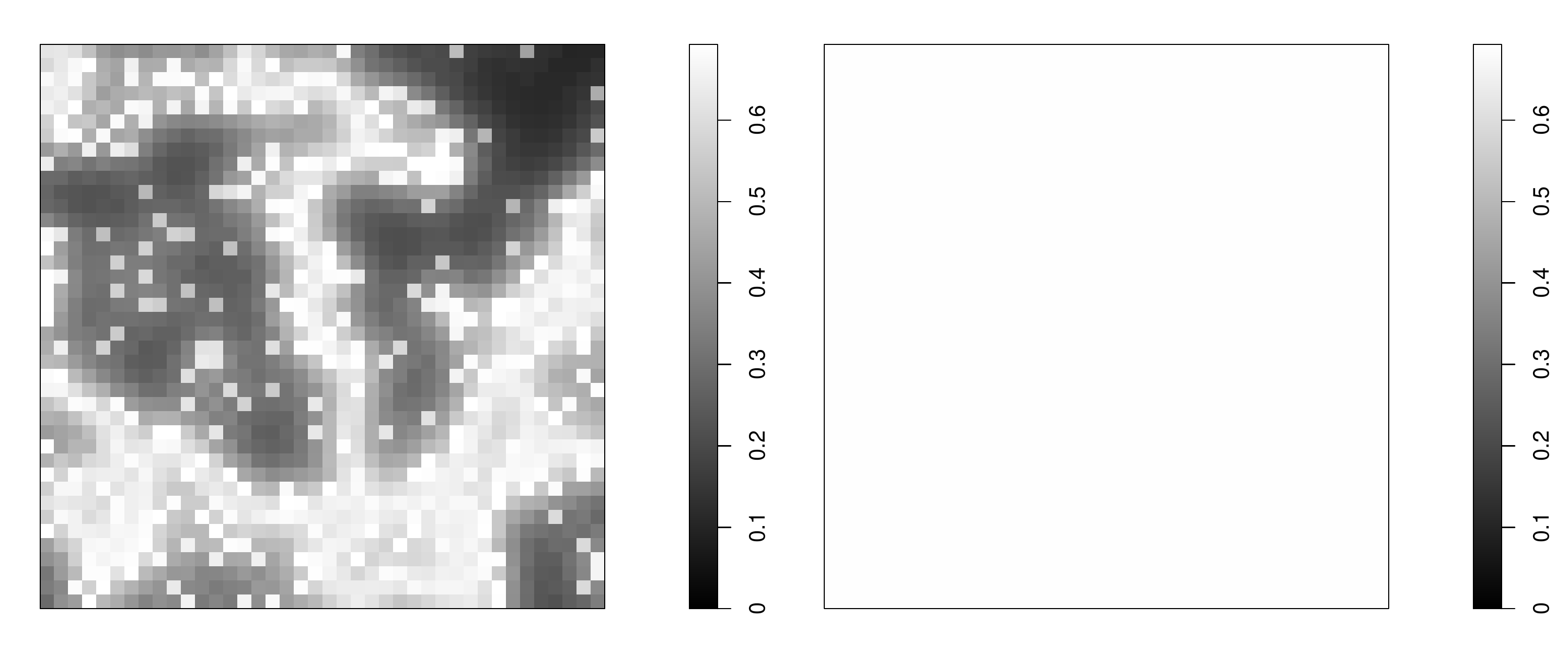}
     \caption{Example of estimated entropy surface for the two scenarios.}\label{fig:entropy}
\end{figure}

Obtaining the entropy surface proceeds as follows. First, the posterior distribution for each $p_u$ is synthesized with the posterior mean. This way, for each scenario and replicate we obtain a single number for $\hat{p}_u$ on every cell. Then, an entropy value is computed over each cell following Equation (\ref{eq:ent_est}) and a smooth spatial function is produced. An example shown in Figure \ref{fig:entropy}, where values range from $0$ (dark areas in the figure) to $\log(2)$ (white areas). The clustered situation (left panel) shows a smoothly varying surface. By comparing the left panels of Figure 1 and \ref{fig:entropy}, one can see that the entropy surface takes low values in areas where pixels are of the same type: white pixels in the top-right part in Figure 1, and black pixels in the bottom-right part of Figure 1, correspond to the darker areas of Figure \ref{fig:entropy} where entropy values are low. In the areas where white and black pixels mix, the entropy surface tends to higher values (whiter areas in Figure \ref{fig:entropy}). The random configuration (right panel of Figure \ref{fig:entropy}) has a constant entropy close to the maximum $\log(2)$; this is expected, as there is no spatial correlation influencing the entropy surface in this scenario. Therefore, such spatial functions properly estimates the entropy of the underlying spatial process. Thanks to the availablity of the marginal posterior distribution of all parameters, any other useful information or synthesis is straightforward to compute. 

\section{Concluding remarks}
\label{sec:disc}

In this paper, we describe an approach to entropy estimation which starts from rigorous posterior evaluation of its components, i.e. the probabilities. This way, we frame entropy within the theory of Bayesian models for spatial data, thus assembling all the available results in this field.

Results from the simulation study enforce the validity of the approach in providing good estimates for the distribution parameters and, consequently, for entropy. The flexibility of the Bayesian paradigm allows to synthesize the posterior distribution of entropy as wished, in order to answer different potential questions. 

Our procedure ensures realistic results, since, when the behaviour of a spatial process is under study, the basic hypothesis is that it is not constant but smoothly varying over space. In the same spirit, an appropriate spatial entropy measure is not a single number, rather it has to be allowed to vary over space as a smooth function.

\bibliographystyle{chicago}
\bibliography{bibdatabase_entropy2}
\end{document}